\documentclass{hjm1}
\usepackage{verbatim}
\numberwithin{equation}{section}
\journalvol{32}
\journalno{1}
\journalyear{2006}

\usepackage{amssymb} 
\usepackage{amssymb} 
\usepackage{latexsym}
\usepackage{amsmath} 
\usepackage{amscd}

\newcommand{\R}{\mbox{$\Bbb R$}}
\def\n{\nabla}

\def\pr#1{{\Bbb P}^{#1}}

\def\X#1{\mathfrak{X}(#1)}

\def\h{\mathfrak{h}}

\def\p#1{\frac{\partial}{\partial #1}}

\def\Pe{P_e}

\def\Pc{P_c}
\def\Pg{P_g}
\def\Pl{P_L}
\def\Pf{{P_{\scriptscriptstyle F}}}

\def\Pfff{{P_{\scriptscriptstyle F}^2}}
\def\Pl{P_{\scriptscriptstyle L}}

\def\PH{P_{\mathfrak{h}}}
\def\pr{\mathfrak{h}}

\def\TM{{\mathcal T}M}

\def\P{{{\mathcal R}}}
\def\H{\mathcal{H}}
\def\h{\mathfrak{h}}

\def\l{\ell}
\def\L{\mathfrak{L}}

\def\lfleche{\smash{\mathop{\vbox{\hbox to 8mm {{\rightarrowfill}}}}
    \limits^{}_{}}}
\def\nn{\noindent}

\newtheorem{idfn}{Remark}

\newtheorem{theorem}{Theorem}
\newtheorem{proposition}{Proposition}
\newtheorem{definition}{Definition}
\newtheorem{lemma}{Lemma}

\newenvironment{remark_num}{\begin{idfn}
\rm}{\end{idfn}}

\begin{document}
\title[Finsler metrizability]{The Euler-Lagrange PDE and Finsler
  metrizability} 

\author[Z.~Muzsnay]{Zolt\'an Muzsnay} 

\address[]{Institute of Mathematics, University of Debrecen, Debrecen,
  H-4010, PBox 12, Hungary} \email{muzsnay@math.klte.hu}




\begin{abstract}
  In this paper we investigate the following question: under what
  conditions can a second-order homogeneous ordinary differential
  equation (spray) be the geodesic equation of a Finsler space.  We
  show that the Euler-Lagrange partial differential system on the
  energy function can be reduced to a first order system on this same
  function. In this way we are able to give effective necessary and
  sufficient conditions for the local existence of a such Finsler
  metric in terms of the holonomy algebra generated by horizontal
  vector-fields.  We also consider the Landsberg metrizability problem
  and prove similar results.  This reduction is a significant step in
  solving the problem whether or not there exists a non-Berwald
  Landsberg space.
\end{abstract}
 
\maketitle

\section{Introduction}
\label{sec:introduction}

A Finsler structure on an $n$-manifold $M$ is a nonnegative function
$F:TM \to \R$ that is smooth and positive away from the zero section
of $TM$, positively homogeneous of degree 1, and strictly convex on
each tangent space.  The energy function $E:TM \to \R$ associated to a
Finsler structure $F$ is defined as $E:=\frac{1}{2}F^2$. This is a
direct generalization of a Riemannian structure. The fundamental
tensor $g_{_E}$ associated to $E$ is formally analogous to the metric
tensor in Riemannian geometry.  It is defined by
\begin{equation} 
  \label{eq:g}
  (g_{_E})_{ij} := \frac{\partial^2 E}{ \partial {y^i} \partial
    {y^j}},
\end{equation} 
in an induced standard coordinate system $(x,y)$ on $TM$.

As in Riemannian geometry, a canonical connection $\Gamma$ can be
defined for a Finsler space \cite{Gri}.

However, since the energy function is not necessarily quadratic and
only homogeneous, the connection is in general non-linear.  We mention
two special types of Finsler spaces: \textit{Berwald spaces}, where
the connection $\Gamma$ is linear, and \textit{Landsberg spaces},
where the connection $\Gamma$ is metric, i.e.~the parallel transport
preserves the norm defined by $g_{\scriptscriptstyle E}$.

Suppose that $M$ is an $n$-manifold endowed with a Finsler structure.
The geodesics are the extremals of the variational problem in which
the Lagrangian is the energy function.  Since $g_{\scriptscriptstyle
  E}$ is non-degenerate, the parametrization of the extremals is
fixed.  The geodesic equation associated to a Finsler structure is
described by the Euler-Lagrange equations
\begin{equation}
  \label{eq:E-L_loc}
  \hphantom{\qquad i=1,...,n}
  \frac{d}{dt}\frac{\partial E}{\partial y^i }- \frac{\partial E}{\partial
  x^i} = 0, \qquad i=1,...,n. 
\end{equation}
Recently several papers were devoted to the problem of characterizing
second-order differential equations coming from a Finsler, a special
Finsler, or a generalized Finsler structure (see for example
\cite{Bao_Shen}, \cite{bryant}, \cite{Kozma_1}, \cite{ok}, \cite{CST},
\cite{shen}, \cite{szil_vat}).  In this paper we offer a contribution
to the solution of this problem.  Now we formulate the problem from
our point of view.

\begin{definition}
  A second-order differential equation on $M$, locally given by 
  \begin{equation}
    \label{eq:sode}
    \hphantom{\qquad i=1,...,n} \ddot{x}^i=f^i(x, \dot{x}), \qquad
    i=1,...,n,
  \end{equation}
  where the functions $f^i$ are positive homogeneous of degree 2 in
  the $\dot{x}$ variable, is called Finsler metrizable, if there
  exists a Finsler structure whose geodesics are described by
  (\ref{eq:sode}).  Moreover, (\ref{eq:sode}) is Landsberg metrizable,
  if it is Finsler metrizable, and in addition we also have
  \begin{equation}
    \label{eq:E_h_metric_loc}
    \frac{\partial g_{jk}}{\partial x^i } - \Gamma_i^l \frac{\partial
      g_{jk}}{\partial y^l} - \Gamma_{ik}^l g_{l j} - \Gamma_{ij}^{l}
    g_{lk}= 0,
  \end{equation}
  where $\Gamma^i_j:= -\frac{1}{2} \frac{\partial f^i}{\partial y^j}$
  are the components of the connection $\Gamma$ associated to
  (\ref{eq:sode}), $\Gamma^i_{jk}:= \frac{\partial
    \Gamma^i_j}{\partial y^k}$, and $g_{ij}=(g_{_E})_{ij}$.
\end{definition}

It follows that a second-order system (\ref{eq:sode}) is Finsler
metrizable if and only if there exists a function $E:TM \to \R$
(energy function), so that
\begin{enumerate}
\item $E$ is homogeneous of degree 2,
\item $E$ is a solution of the Euler-Lagrange system
  (\ref{eq:E-L_loc}) considered as a second-order partial differential
  equation with respect to $E$,
\item the quadratic form $g_{\scriptscriptstyle E}$ defined by
  (\ref{eq:g}) is positive definite.
\end{enumerate}

Euler's theorem for homogeneous functions implies that the homogeneity
condition on $E$ can be described by the equation
\begin{equation}
  \label{eq:E_hom_loc}
  y^i \frac{\partial E}{\partial y^i}-2E = 0.
\end{equation}
The Euler-Lagrange partial differential equation associated to
(\ref{eq:sode}) is
\begin{equation}
  \label{eq:E_EL_loc}
  \hphantom{\qquad j=1,..,n} y^j \frac{\partial^2 E}{\partial x^j
    \partial y^i}+f^j \frac{\partial^2 E}{\partial y^j \partial y^i}-
  \frac{\partial E}{\partial x^i } = 0, \qquad i=1,..,n,
\end{equation}
in induced local coordinates $(x,y)$ on $TM$. We arrive at the
reformulation of the metrizability property in terms of a partial
differential system:

\begin{proposition}
  \label{theo:reformulation}
  A second-order differential equation (\ref{eq:sode}) is
  \begin{enumerate}
  \item \textit{Finsler metrizable}, if and only if, there exists a
    solution $E:TM \to \R$ to the second-order PDE system formed by
    the equations (\ref{eq:E_hom_loc}) and (\ref{eq:E_EL_loc}) so that
    the quadratic form $g_E$ defined in (\ref{eq:g}) is positive
    definite;
  \item \textit{Landsberg metrizable}, if and only if there exists a
    solution $E:TM \to \R$ to the third-order\footnote{The equation
      (\ref{eq:E_h_metric_loc}) is a 3rd order PDE, taking into
      account of (\ref{eq:g}).} PDE system (\ref{eq:E_h_metric_loc}),
    (\ref{eq:E_hom_loc}) and (\ref{eq:E_EL_loc}) so that the quadratic
    form $g_E$ is positive definite.
  \end{enumerate}
\end{proposition}

\noindent
The main results of this paper can be found in Sections
\ref{cha:finsler} and \ref{cha:landsberg}.

\smallskip

In Section \ref{cha:finsler} we consider the problem of Finsler
metrizability.  Using the integrability conditions of the
corresponding PDE, we show that the system is equivalent to a first
order PDE on the same unknown function (Theorem \ref{theo:finsler}).
We formulate a necessary and sufficient condition for the local
metrizability in terms of a distribution $\H$ associated to the spray
(Theorem \ref{theo:C_D} and \ref{theo:finsler_2}).  $\H$ is called
\textit{holonomy distribution} or \textit{holonomy algebra}
\cite{Kle}, and it is generated by the horizontal vector fields and
their successive Lie-brackets. 

In Section \ref{cha:landsberg} we consider the problem of Landsberg
metrizability.  We show that the corresponding third-order system can
be reduced to a first order PDE on the same energy function (Theorem
\ref{theo:landsberg}).  As in the previous case, we are able to
formulate a necessary and sufficient condition for the metrizability
in terms of a distribution $\L$ (Theorem \ref{theo:landsberg_2}).  The
distribution $\L$ is generated by the holonomy algebra and the image
of the Berwald curvature.

\medskip

In Sections \ref{sec:remarks} and \ref{sec:exist-non-berw} we
illustrate some consequences of the results on Finsler and Landsberg
metrizability.  We also discuss the famous problem of whether there
exists a non-Berwald Landsberg space.  As we show through several
examples, Theorem \ref{theo:landsberg_2} offers a promising
alternative approach to solve this problem.

\section{Preliminaries}

\subsection{Notations, conventions}

Throughout this paper $M$ will denote an $n$-dimensional smooth
manifold. $C^\infty(M)$ denotes the ring of real-valued smooth functions,
$\mathfrak X(M)$ is the $C^\infty(M)$-module of vector fields on $M$,
$\pi:TM \to M$ is the tangent bundle of $M$, $\TM=TM \setminus 0$ is
the slit tangent space.  We will essentially work on the manifold
$TM$ and on its tangent space $TTM$.  When there is no danger of
confusion, $TTM$ and $T^*TM$ will simply be denoted by $T$ and $T^*$,
respectively.  $T^v=\mathrm{Ker}\,\pi_*$ will be the vertical sub-bundle
of $T$.

The exterior differential, the Lie differential (with respect to $X
\in \mathfrak X (M)$) and the interior product (induced by $X$) are
denoted by $d$, ${\mathcal L}_X$ and $i_X$, respectively.

We denote by $\Lambda^k(M)$ and $S^k(M)$ the $C^\infty (M)$-modules of
the skew-symmetric and symmetric $k$-forms.  The Fr\"olicher-Nijenhuis
theory provides a complete description of the derivation of
$\Lambda(M)$ with the help of vector-valued differential forms, for
details we refer to \cite{FN}.  The $i_*$ and the $d_*$ type
derivation associated to a vector valued $l$-form $L$ will be denoted
by $i_L$ and $d_L$. They can be defined in the following way:
\begin{enumerate}
\item if ${\rm deg}\,L = 0$, i.e. $L \in \mathfrak{X}(M)$, then $i_L
  \omega : = \omega(L)$, \ and $d_L\omega : = {\mathcal L}_L\omega$;
  
\item if ${\rm deg} \ L =l>1$, then
   
  \quad
  \begin{math}
    i_L\omega(X_1,...,X_l) : =  \omega(L(X_1,...,X_l)),
  \end{math}
  \quad{}   \ \ \ \ for $\omega \in \Lambda^1(M)$;

  \quad 
  \begin{math}
    d_Lf(X_1,\cdots,X_l) : =  df(L(X_1,\cdots,X_l)),
  \end{math}
  \quad  for $f \in C^\infty(M)$.
\end{enumerate}

\subsection{Geometry associated to a spray}

Let $J$ be the canonical vertical endomorphism of $T~(=TTM)$ and $C \in
\mathfrak X(TM)$ the canonical vertical vector field.  In an induced
local coordinate system $(x^i,y^i)$ on $TM$ we have
\begin{displaymath}
  J = dx^i \otimes \frac {\partial}{\partial y^i },
  \qquad C = y^i \frac {\partial}{\partial y^i }.
\end{displaymath}

\begin{remark_num}
  \label{theo:Pc}
  Using the canonical vector-field, equation (\ref{eq:E_hom_loc})
  can be written in the form $\Pc E = 0$, where
  \begin{math}
    \Pc : C^\infty(TM) \to C^\infty(TM)
  \end{math}   
  is a first-order differential operator defined on a function $E:TM
  \to \R$ by
  \begin{equation}
    \label{eq:E_hom}
    \Pc E:= {\mathcal L}_C E - 2E.
  \end{equation}
\end{remark_num}

A \textit{spray} is a vector field $S \in \mathfrak X(TM)$ on $TM$
satisfying the relations $JS = C$ and $[C,S]=S$. The coordinate
representation of a spray $S$ takes the form
\begin{equation}
  \label{eq:S}
  S = y ^i \frac {\partial}{\partial x^i } +f^i (x,y)
  \frac {\partial}{\partial y^i },
\end{equation}
where $f^i(x,y)$ is positive-homogeneous of degree 2 in $y=(y^j)$. The
integral curves of a spray are curves $\gamma : I \to M$ so that $S
\circ \dot \gamma = \ddot \gamma$.  They are the solutions of the
equations 
\begin{math}
  \ddot{x}^i = f^i(x,\dot{x}).
\end{math}

To every spray $S$ a \textit{connection} $\Gamma: = [J,S]$ can be
associated \cite{Gri}. We have $\Gamma ^2 = \textrm{id}_T$, and the
eigenspace of $\Gamma$ corresponding to the eigenvalue $-1$ is the
vertical space $T^v$. We denote the eigenspace belonging to the
eigenvalue +1 of $\Gamma$ by $T^h$ and we call it the
\textit{horizontal space}. Then 
\begin{displaymath}
  T = T^h \oplus T^v.
\end{displaymath}
The  horizontal and the vertical projector belonging to $\Gamma$ are 
\begin{math}
  h : = \frac{1}{2}(\textrm{I} + \Gamma ),
\end{math}     
and $v := \mathrm{id}_T-h$.  The almost complex structure associated
to $\Gamma$ is the vector valued 1-form $F$ on $TM$ such that $FJ = h$
and $Fh = -J$.  The \textit{curvature} of the connection $\Gamma$ is
the vector-valued 2-form
\begin{equation}
  \label{def:R}
  R : =  -\frac{1}{2}[h,h].
\end{equation}
A linear connection on $TM$, called the {\it Berwald connection}, can
also be associated to $S$. It is defined by:
\begin{displaymath}
  \n \Gamma = 0, \qquad \n_{hX} JY = [ h, JY]X, \qquad \n_{JX} JY = [ J,
  JY]X;
\end{displaymath}
$X, Y \in \mathfrak X (TM)$. In an induced coordinate system $(x,y)$
we have 
\begin{equation}
  \label{berwald}
  \left\{ \quad
    \begin{aligned}{}
      \n_{_{\p{y^i}}} {\scriptstyle \p{y^j}}&=0,
      \\
      \n_{_{ \p{x^i}}} {\scriptstyle \p{y^j}}&= \n_{_{\p {y^j}}}
      {\scriptstyle \p{x^i}} = \Gamma_{i j}^k \, {\scriptstyle \p
        {y^k}},
      \\
      \n_{_{\p{x^i}}} {\scriptstyle \p{x^j}} &= \Gamma_{i j}^k
      {\scriptstyle \p {x^k}}
      + \Bigl( \tfrac{\partial \Gamma_j^l} {\partial x^i}
      + \Gamma_j^k \Gamma_{i k}^l
      - \Gamma_k^l \Gamma_{i j}^k \Bigl)\, {\scriptstyle \p {y^l}}.
    \end{aligned}
      \right.
\end{equation}
where $\Gamma^k_i:= -\frac{1}{2}\frac{\partial f^k}{\partial y^i}$ and
$\Gamma^k_{ij}:= \frac{\partial \Gamma^k_i}{\partial y^j}$.
Considering the $(h,v,v)$ components of the classical curvature of the
Berwald connection we obtain a tensor-field
\begin{equation}
  \label{eq:P}
  \P (X,Y,Z)=\n_{hX}\n_{JY}JZ - \n_{JY}\n_{hX}JZ - \n_{[hX, JY]}JZ
\end{equation}
called the \textit{Berwald curvature} in Shen's monograph \cite{shen}.

\begin{remark_num}
  Using the coordinate expressions (\ref{berwald}), it is easy to see
  that locally we have
  \begin{displaymath}
    \P = -\frac{1}{2} \frac{\partial^3 f^l}{\partial y^i \partial y^j
      \partial y^k} \, dx^i \otimes dx^j \otimes dx^k \otimes \p
    {y^l}.
  \end{displaymath}
  Therefore the connection $\Gamma$ is linear, and the corresponding
  Finsler space is of Berwald type, if and only if, $\P=0$.
\end{remark_num}

\medskip

\begin{remark_num}
  \label{theo:Pg}
  Using the Berwald connection, we can introduce a third-order
  differential operator 
  \begin{math}
    \Pg : C^\infty(TM) \longrightarrow Sec \, (T^* \otimes S^2 T^*),
  \end{math} 
  given by
  \begin{equation}
    \label{eq:E_h_metric_P}
    (\Pg E)(X,Y,Z):=\n_{hX} g_{_E}(JY, JZ),
  \end{equation}
  for $X,Y,Z \in \mathfrak X (TM)$. Then (\ref{eq:E_h_metric_loc})
  takes the form
  \begin{math}
    \label{eq:E_h_metric}
    \Pg E = 0.
  \end{math}
\end{remark_num}

\subsection{Lagrangian and spray}

A Lagrangian $E: TM \to \R$ is called \textit{regular}, if the 2-form
\begin{displaymath}
  \Omega_E : = dd_JE
\end{displaymath}
is symplectic. This holds if and only if
\begin{math}
  \det \Bigl( \frac {\partial^2E}{\partial y^\alpha \partial y^\beta }
  \Bigl) \neq 0.
\end{math}
Let $S \in \mathfrak X (TM)$ be a spray. We introduce a second-order
differential operator
\begin{math}
  \Pe : C^\infty(TM) \to Sec~T^*, 
\end{math} 
given by
\begin{equation}
  \label{eq:E_EL_P}
  \Pe E := i_S \Omega_E +d{\mathcal L}_C E -dE.
\end{equation}
It is not difficult to see that $\Pe E$ is a semi-basic 1-form for all
$E \in C^\infty(TM)$, and its coordinate representation takes the form
$\Pe E = \omega_i \, dx^i$ where the coefficients $\omega_i$ are the
functions appearing in the left-hand side of the Euler-Lagrange
equation (\ref{eq:E_EL_loc}). Therefore $S$ corresponds to the
geodesic equation of $E$ if and only if the equation $\Pe E = 0$
is valid. So we have the

\begin{remark_num}
  \label{theo:Pe}
  If $S$ is a spray, then $\Pe E=0$ is the coordinate-free
  ex\-pres\-sion of the Euler-Lagrange partial differential equation
  (\ref{eq:E_EL_loc}) associated to $S$.
\end{remark_num}

\subsection{Formal integrability}

In order to solve the metrizability problems formulated above, we have
to deal with partial differential systems. We shall use Spencer's
technique of formal integrability in the form explained in
\cite{grif_muzs_2}; for a detailed account see \cite{BCGGG}.  We
recall here only some basic notions in order to fix the terminology.

Let $B$ be a vector bundle over $M$.  If $s$ is a section of $B$, then
$j_{k,p}s=(j_k s)_p$ will denote the $k$th order jet of $s$ at the
point $p \in M$.  The bundle of $k$th order jets of the sections of
$B$ is denoted by $J_kB$.  In particular $J_k(\R_M)$ will denote the
$k$th order jet of the sections of the trivial line bundle, i.e.~the
real valued functions.  If $B_1$ and $B_2$ are two vector bundles over
the same manifold $M$ and
\begin{displaymath}
  P:Sec\,(B_1) \to Sec \,(B_2)
\end{displaymath}
is a linear differential operator of order $k$, then the morphism
$p_{k+l}(P):\, J_{k+l}(B_1) \to J_l(B_2)$ defined by
\begin{displaymath}
  \hphantom{\qquad l=0,1,2,...} p_{k+l}(P) \,
  \bigl(j_{k+l,p}(s)\bigl):=j_{l,p}(Ps), \qquad l=0,1,2,...
\end{displaymath}
is called the $l$th order prolongation of $P$.  $R_{k+l,
  p}(P):=\mathrm{Ker}\, p_{k+l}(P)_p$ will denote the bundle of the
formal solutions of order $k+l$ at $p$.  A differential operator $P$
is called \textit{formally integrable} at $p\in M$, if $R_{k+l}(P)$ is
a vector bundle for all $l \geq 0$, and $\overline{\pi}_{k+l,p} :
R_{{k+l},p}(P)\rightarrow R_{k+l-1,p}(P)$ is onto for every $l\geq 1$.
In analytical terms, formal integrability implies for arbitrary
initial data the existence of solutions (see.~\cite{BCGGG}, p.~397).

\smallskip

$\sigma_{k}(P):S^{k}T^*M \otimes B_1 \to B_2$ is the symbol of $P$,
defined as the highest order terms of the operator, and
$\sigma_{k+l}(P):S^{k+l}T^*M \otimes B_1 \to S^lT^*M \otimes B_2$ is
the symbol of the $l$-th order prolongation of $P$.  We write
\begin{alignat*}{2}
  & g_{k,p} (P) & & = \mathrm{Ker} \, \sigma_{k,p}(P),
  \\
  & g_{k,p} (P)_{e_1...e_j} & & = \bigl\{A\in g_{k,p}(P) \mid i_{e_1}A
  = ....  = i_{e_j}A = 0 \bigl\}, \quad j=1,...,n,
\end{alignat*}
where $\{e_1,...,e_n\}$ is a basis of $T_p M$.  A basis
$\{e_i\}_{i=1}^n$ of $T_pM$ is called \textit{quasi-regular} if
\begin{displaymath}
  \textrm{dim} \, g_{k+1,p}(P) = \textrm{dim} \, g_{k,p}(P) + {\sum
    _{j=1}^{n}} \textrm{dim} \, g_{k,p}(P)_{e_1...e_j} .
\end{displaymath}
A symbol is called \textit{involutive}\footnote{There is a slight
  problem of language here. In the works of Cartan, and more generally
  in the theory of exterior differential systems, ``involutivity''
  means more than the existence of a quasi-regular basis and it refers
  to "integrability" (cf.  \cite{BCGGG}, p.107, 140). Here we are
  following the terminology of Goldschmidt (cf.  \cite{BCGGG},
  p.\,409).}  at $p$, if there exists a quasi-regular basis at $p$.
The notion of involutivity allows us to check the formal integrability
in quite a simple way:

\bigskip

\nn \textbf{Theorem}
\label{theo:cartan_kahler}  (Cartan-K\"ahler). \textit{Let $P$ be a linear
  partial differential o\-pe\-ra\-tor.  Suppose that $g_{k+1}(P)$ is
  regular, i.e.~$R_{k+1}(P)$ is a vector bundle on $R_k(P)$.  If the
  map $\overline {\pi }_{k} : R_{k+1}(P)\longrightarrow R_k(P)$ is
  onto and the symbol is involutive, then $P$ is formally integrable.}

\bigskip

\section{Finsler metrics with prescribed geodesics}
\label{cha:finsler} 

\medskip

In this paragraph we are going to investigate the following problem:
\textit{under which conditions can a second order differential
  equation (\ref{eq:sode}) be the geodesic equation of a Finsler
  metric}.  As we explained in Section \ref{sec:introduction}
(Proposition \ref{theo:reformulation}) we have to look for a solution
of the PDE comprised of (\ref{eq:E_hom_loc}) and (\ref{eq:E_EL_loc}).
Therefore we have to deal with the second-order system
\begin{equation}
  \label{eq:e_L_hom}
  \Pf := (\Pc , \, \Pe )
\end{equation}
where $\Pc$ and $\Pe$ are defined in (\ref{eq:E_hom}) and
(\ref{eq:E_EL_P}).  We will prove the following theorems:

\begin{theorem}
  \label{theo:finsler}(Reduction of $\Pf$.)
  A Lagrangian $E: TM \to \R$ is a solution of the second order
  operator $\Pf$, if and only if, it is a solution of the first order
  system
  \begin{equation}
    \label{eq:syst_2_1}
    \left\{\quad 
      \begin{aligned}
        {\mathcal L}_CE-2E&=0,
        \\
        d_\h E & =0,
        \end{aligned}
    \right.
  \end{equation}
  where ${\mathcal H} \subset T (=TTM)$ is the holonomy algebra
  generated by the horizontal vector fields and their successive
  Lie-brackets, and $\h:T\to \H$ is an arbitrary projection on $\H$.
\end{theorem}

\noindent
\begin{idfn}
  For $X \in \X {TM}$ we have
  \begin{math}
    d_\h E (X) = \h X (E) = {\mathcal L}_{\h X} E,
  \end{math}
  so the second equation of (\ref{eq:syst_2_1}) means simply that the
  Lie-derivative of $E$ with respect to vector-fields in the holonomy
  distribution $\H= \mathrm{Im}\, \h$ is zero.  This property is
  independent of the projection $\h$ of $\H$ chosen.  
\end{idfn}

\textit{Proof of Theorem \ref{theo:finsler}.}
Let us suppose that $E:TM \to \R$ is a solution of
(\ref{eq:syst_2_1}). Since $T^h \subset \H$, we have $\h \circ h = h
$. Therefore
\begin{displaymath}
  d_hE  = d_{\h \circ h}E = i_h d_\h E  - d_\h  i_h E + i_{[h, \h]} E
  =   i_h d_\h E = 0
\end{displaymath}
since the action of an $i_*$-type derivation is trivial on functions.
Moreover as $S$ is homogeneous, $hS=S$ and
\begin{displaymath}
  {\mathcal L}_{S}E= {\mathcal L}_{hS}E=d_h E (S) =0.
\end{displaymath}
Writing the Euler-Lagrange operator in the form
\begin{displaymath}
  \Pe E = i_Sdd_JE + d{\mathcal L}_CE - dE = d_J{\mathcal L}_SE -i_{[J,S]}dE =
  d_J{\mathcal L}_SE -2d_hE
\end{displaymath}
we obtain that $\Pe E =0$ and $E$ is a solution of (\ref{eq:e_L_hom}).

Let us suppose now that $E:TM \to \R$ is a solution of
(\ref{eq:e_L_hom}). We have
\begin{equation}
  \label{eq:el}
  i_S \Omega_E = d(E - {\mathcal L}_CE) = -dE.
\end{equation}
Since $[J,J] = 0$, we have
\begin{math}
  d^2_J = d_J \circ d_J = d_{[J,J]}=0,
\end{math}
and 
\begin{math}
  i_J \Omega_E = 0,
\end{math}
 so
\begin{equation}
  \label{eq:ic}
  i_C \Omega_E=i_{JS} \Omega_E = i_Si_J\Omega_E - i_J i_S \Omega_E =
  i_J dE.
\end{equation}
On the other hand, for every $X \in \X{TM}$ we have 
\begin{displaymath}
  i_S \Omega_E(X) = \Omega_E(S, X) = - \Omega_E(C, FX) = - i_F i_C
  \Omega_E (X),
\end{displaymath}
i.e.
\begin{equation}
  \label{eq:is}
  i_S \Omega_E = i_F i_C \Omega_E .
\end{equation}
Putting (\ref{eq:ic}) into (\ref{eq:is}) we obtain
\begin{equation}
  \label{eq:elh}
  i_S \Omega_E = - i_Fi_C \Omega_E = -i_F i_J dE = -d_v E = -d E + d_h
  E.
\end{equation}
Comparing (\ref{eq:elh}) with (\ref{eq:el}) we obtain that 
\begin{math}
  d_h E=0.
\end{math}
It follows that $hX(E)=0$, i.e.~$E$ is constant with respect
to horizontal vector fields. Therefore it must be constant on the
distribution generated by the horizontal sub-bundle taking the
recursive Lie-bracket operations, i.e.~on $\H$.  This means that we
have $d_\h E =0$ and $E$ is a solution of (\ref{eq:syst_2_1}).
\hfill{$\Box$}

\begin{remark_num}
  \label{rem:num_1}
    $E$ is a solution of (\ref{eq:syst_2_1}) if and only if it is
    a solution of
    \begin{equation}
      \tag{\ref{eq:syst_2_1}'}
      \label{eq:syst_2_1_1}
      \left\{\quad
        \begin{aligned}
          {\mathcal L}_CE-2E&=0,
          \\
          d_h E & =0,
        \end{aligned}
      \right.
    \end{equation}
    where $h$ is simply the horizontal projection associated to
    $\Gamma$, so (\ref{eq:syst_2_1}) and (\ref{eq:syst_2_1_1}) are
    equivalent. However, as we will see in Proposition
    \ref{theo:D_int}, under regularity assumption the system
    (\ref{eq:syst_2_1}) is integrable while (\ref{eq:syst_2_1_1}) is
    not, unless the curvature is zero.  Indeed, we have
    \begin{displaymath}
      d_R E= - \tfrac{1}{2} d_{[h,h]}E = - \tfrac{1}{2} d_h d_h E,
    \end{displaymath}
    therefore $d_R E=0$ is a compatibility condition for
    (\ref{eq:syst_2_1_1}).
\end{remark_num}
    
\begin{remark_num}
  \label{rem:num_2}
   Let us introduce the first order differential operator $\PH
    :C^\infty(TM) \longrightarrow Sec\,(T^*)$ by the rule
    \begin{equation}
      \label{eq:d_p}
      \PH E\, (X):=\pr X(E), 
    \end{equation}
    $E \in C^\infty(TM)$, $X\in \mathfrak{X}(TM)$, and the differential
    operator 
    \begin{equation}
      \label{eq:D_3}
      \Pfff:=(\Pc, \, \PH)
    \end{equation}
    corresponding to the system (\ref{eq:syst_2_1}).  Theorem
    \ref{theo:finsler} shows that a Lagrangian is a solution of $\Pf$
    if and only if it is a solution of $\Pfff$.
\end{remark_num}

\begin{theorem}
  \label{theo:C_D}
  Let $S$ be a spray over the manifold $M$.  If $C \in \H$, then there
  is no Finsler metric whose geodesics are given by $S$.
\end{theorem}

\begin{proof} 
  Let $S$ be a spray and $E:TM \to \R$ a Lagrangian.  From Proposition
  \ref{theo:reformulation} we know that if $E$ is an energy function
  associated to $S$, then it is a solution of $\Pf=(\Pc, \, \Pe)$, and
  by Theorem \ref{theo:finsler} we obtain that $E$ satisfies the
  equations $\mathcal{L}_CE -2E = 0$ and $d_\h E =0$.  If $C\in \H$,
  we have also
  \begin{displaymath}
    0=\PH E (C)=(\h C) E = CE ={\mathcal L}_CE,
  \end{displaymath}
  therefore $E = 0$. Since $E$ has to be a regular Lagrangian, this is
  impossible and the proposition is proved.
\end{proof}

Let us consider the case when $C \not \in \H$.  We have the following 
\begin{theorem}
  \label{theo:finsler_2}
  Let $S$ be an analytical spray over the analytical manifold $M$.  If
  $C \not \in \H$ and $\H$ has constant rank in a neighbourhood of
  $v\in \TM$, then there exists an analytical Finsler metric in a
  neighbourhood of $v$ such that the geodesics are given by $S$ if and
  only if the kernel of the first prolongation of (\ref{eq:syst_2_1})
  at $v$ contains positive definite initial data.
\end{theorem}

\begin{remark_num}
  Let $(x^i)$ be a local coordinate system on $M$, $(x^i,y^i)$ the
  associated coordinate system on $TM$ in the neighborhood of $v$.  If
  $p:=j_k(E)_v \in J_{2}(\R_{TM})$ is a $k$th order jet of a real
  valued function $E$ on $TM$ we set
  \begin{equation}
    \label{eq:ch5_1}
    s_{i_1...i_a \underline{i_{a+1}..i_l}}(p):= \frac{\partial^l E}{\partial
      x^{i_1} ... \,  \partial x^{i_a} \partial y^{i_{a+1}} ... \,   \partial
      y^{i_l}}(v), \qquad 1 \leq l \leq k. 
  \end{equation}
  Then
  \begin{math}
    (s,s_j,s_{\underline{j}}, s_{jk}, s_{j\underline{k}},
    s_{\underline{jk}})
  \end{math}
  gives a coordinate system on $J_{2,v}(\R_{TM})$.
  Using the notation of (\ref{eq:D_3}) introduced in Remark
  \ref{rem:num_2}, positive definite initial data for the first
  prolongation of (\ref{eq:syst_2_1}) at $v$ is simply an element
  $s_{2,v} \in J_{2,v}(\R_{TM})$ represented as
  \begin{math}
    s_{2,v}=(s,s_i,s_{\underline{i}}, s_{ij}, s_{i\underline{j}},
    s_{\underline{ij}}) \in \R^{1+ (n+n)+ \frac{n(n+1)}{2}+
      n^2+\frac{n(n+1)}{2}}
  \end{math}
  such that ${(s_{\underline{ij}})}_{1 \leq i,j\leq n}$ determines a
  positive definite quadratic form, and $s_{2,v}$ is a second order
  solution of $\Pfff$ at $v$. This last condition gives \textit{linear
    algebraic} equations on the coordinates of $s_{2,v}$.
\end{remark_num}

\textit{Proof of Theorem \ref{theo:finsler_2}.} The proof is based on
Theorem \ref{theo:finsler} and on Proposition \ref{theo:D_int} proved
below.  Indeed, if $\Pfff$ is formally integrable (see Proposition
\ref{theo:D_int}), then for every initial condition we have an
infinite order formal solution of $\Pfff$. In the analytic case, this
formal solution gives an analytical solution in an open neighborhood
of $\TM$.  Theorem \ref{theo:finsler} shows that this solution is also
a solution of the operator $\Pf$. In this way we obtain an analytical
solution of $\Pf$, i.e.~a homogeneous function which satisfies the
Euler-Lagrange equation associated to $S$.  Therefore $S$ is locally
Finsler metrizable.  \hfill{$\Box$}

\begin{proposition}
  \label{theo:D_int}
  Let $S$ be a spray over $M$ so that $C \not \in \H$ and the rank of
  $\H$ is locally constant. Then the differential operator
  $\Pfff=(\Pc, \, \PH)$ is formally integrable.
\end{proposition}

\begin{proof}
  First of all remark that $\Pfff$ is a regular differential operator
  because, by the hypothesis, rank of $\H=\mathrm{Im}\, \h$ is locally
  constant.  Moreover, using Lemma \ref{lemma_1}, Lemma \ref{lemma_2}
  and the Cartan-K\"ahler theorem on formal integrability (see page
  \pageref{theo:cartan_kahler}), we obtain the proposition.
\end{proof}

\begin{lemma}
  \label{lemma_1}
  Every first order solution of $\Pfff$ can be lifted to a second
  order solution.
\end{lemma}

\begin{proof}
  It is easy to see from their local description that the symbol of
  $\Pc$ and $\PH$ can be interpreted as a map
\begin{alignat*}{2}
  & \sigma_1(\Pc) : T^* \to \R, & & \sigma_1(\Pc)B_1 = B_1(C)
  \\
  & \sigma_1(\PH) : T^* \to T^* \qquad & & (\sigma_1(\PH)B_1) (X) =
  B_1(\pr X)
\end{alignat*}
for all $B_1 \in T^*$, $X \in T$. The symbol of the first
prolongations are defined by
\begin{alignat*}{2}
  & \sigma_2(\Pc) : S^2T^* \to T^*, & & (\sigma_2(\Pc)B_2) (X) =
  B_2(X, C)
  \\
  & \sigma_2(\PH) : S^2T^* \to T^* \otimes T^* \qquad & &
  (\sigma_2(\PH)B_2) (X,Y) = B_2(X, \pr Y)
\end{alignat*}
for all $B_2 \in S^2T^*$, $X,Y \in T$.  Comparing the first
prolongation of the symbols, we can easily find that for every $B_2
\in S^2 T^*$ and $X\in T$ we have
  \begin{displaymath}
    (\sigma_2 (\Pc)B_2)(\pr X) - (\sigma_2 (\PH)B_2)(C,X) = B_2(\pr X,
    C ) - B_2(C, \pr X) = 0,
  \end{displaymath} 
  and there is no more relation between the two symbols. That is, if
  we consider the map
  \begin{math}
    \tau: T^* \oplus (T^* \otimes T^*) \longrightarrow T^*
  \end{math}
  defined for $B_1 \in T^*$, $B_2 \in T^* \otimes T^*$ and for $X \in
  T$ as
  \begin{displaymath}
    \bigl(\tau (B_1, B_2)\bigl) (X):= B_1 (\pr X) - B_2 (C, X),
  \end{displaymath}
  then we find the commutative diagram:

  \begin{displaymath}
    \begin{CD}
      & 0 && 0 && &0
      \\
      & \downarrow && \downarrow & & & \downarrow
      \\
      0 \longrightarrow & \ g_2(\Pfff) @>i>> S^2 T^*
      @>\sigma_2(\Pfff)>> & T^* \oplus (T^* \otimes T^*) @>\tau>> T^*
      & \longrightarrow 0
      \\
      & @VVV @VVV & @VVV
      \\
      0 \longrightarrow & R_2 @>i >> J_2(\R_{TM}) @>p_2(\Pfff)>> &
      J_1(\R_{TM} \oplus T^*)
      \\
      & @VV\overline{\pi}V @VV\pi V & @VV\pi V
      \\
      0 \longrightarrow & R_1 @>i >> J_1(\R_{TM}) @>p_1(\Pfff)>> & \R
      \oplus T^*
      \\
      & && @VVV & @VVV
      \\
      & && 0 & & & 0
    \end{CD}
  \end{displaymath}
  where the successive arrows represent exact sequences.  ($R_1$ and
  $R_2$ denote the spaces of the first and second order formal
  solutions of $\Pfff$.)
  
  Every first order solution of $\Pfff$ can be lifted into a
  second-order formal solution if and only if the map $\overline{\pi}:
  {\mathcal R}_2 \to {\mathcal R}_1$ is onto.  We know by a lemma of
  homological algebra that there exists a map
  \begin{math}
    \varphi : R_1 \longrightarrow T^* \, (= \mathrm{Im}\,\tau)
  \end{math}
  such that
  \begin{equation}
    \label{eq:ker_im}
    \mathrm{Im}\, \overline{\pi} = \mathrm{Ker} \,\varphi.
  \end{equation}
  This map can be constructed for a first order formal solution
  $j_{1,v}(E)\in R_1$ of $\Pfff$ at $v \in \TM$ as follows:
  \begin{equation}
    \label{varphi_pi}
    \varphi_v(E) : = \tau (\n \Pfff E)_v.
  \end{equation}
  Let us compute how this map acts. If $E:TM \to \R$ is a function
  such that $j_{1,v}(E) \in R_{1,v}$, then $(\Pfff E)_v=0$, that is
  $(\PH E)_v=0$ and $(\Pc E)_v=0$.  Evaluating $\varphi_v (E)$ on an
  arbitrary vector $X\in T$ we find that
  \begin{displaymath}
    \label{varphi}
    \begin{aligned}
      &\varphi_v (E)(X) = \tau (\n \Pfff E)_v(X) = (\n \Pc E)_v(\h X)-
      (\n \PH E)_v (C, X)
      \\
      & = \bigl({\mathcal L}_{\pr X}({\mathcal L}_CE-2E)-{\mathcal
        L}_C({\mathcal L}_{\pr X} E)\bigl)_v = \bigl({\mathcal L}_{\pr
        X}({\mathcal L}_{C}E) -{\mathcal L}_{C}({\mathcal L}_{\pr X}E)
      -2 {\mathcal L}_{\pr X}E\bigl)_v
      \\
      & = \bigl({\mathcal L}_{[\pr X,C]}E\bigl)_v -2 (\PH E)_v(X) =
      \bigl({\mathcal L}_{[\pr X,C]}E\bigl)_v .
    \end{aligned}
  \end{displaymath}
  Now we can remark, that if $X \in \X {TM}$, then $[\h X, C] \in \H
  $.  Indeed $\H$ can be generated by the successive brackets of the
  horizontal basis
  \begin{math}
    \{h_1, ... h_n\},
  \end{math}
  where
  \begin{math}
    h_i := h \bigl(\p {x^i}\bigl)= \p {x^i} - \Gamma_i^\alpha \p
    {y^\alpha}.
  \end{math}
  Since $S$ is homogeneous, we have $[h_i, C]=0$. By the Jacobi
  identity, this is also true for the successive brackets of the
  $h_i$'s.  If we consider an arbitrary $Y \in \H$, then it can be
  written as a linear combination of the elements $Y=g^\alpha
  Y_\alpha$, where $Y_\alpha$ can be obtained by successive brackets
  of the $h_i$'s.  Thus we have
  \begin{displaymath}
    [Y,C]=[g^\alpha Y_\alpha, C]=-(Cg^\alpha) Y_\alpha +
    g^\alpha[Y_\alpha,C] = -(Cg^\alpha)Y_\alpha,
  \end{displaymath}
  which shows that $[Y,C] \in \H$.
  
  Continuing the above computation of $\varphi_v (E)$ we find that
 \begin{displaymath}
   \varphi_v (E)(X) = \bigl({\mathcal L}_{[\pr X,C]}E\bigl)_v =
   \bigl({\mathcal L}_{\h [\pr X,C]}E\bigl)_v = (\PH E )_v([\pr
   X,C])=0,
  \end{displaymath}
  since $(\PH E)_v$ vanishes on $\H_v$. It follows that $\varphi_v$ is
  identically zero, by (\ref{eq:ker_im}) we conclude that
  $\mathrm{Ker} \, \varphi_v = R_1$, and $\overline{\pi}$ is onto.
  Hence every first order solution of $\Pfff$ can be lifted into a
  second order solution.
\end{proof}
  \medskip{}
  
\begin{lemma}
  \label{lemma_2}
  The symbol of $\Pfff$ is involutive.
\end{lemma}

\begin{proof}
  Let $k$ be the co-dimension of $\H$. Since $\mathrm{dim} \, \H \geq
  n$ and $C \not \in \H$, we have $n \leq \mathrm{dim} \, \H \leq 2n-1
  $ and $1 \leq k \leq n$. Let us consider the basis
  \begin{equation}
    \label{quasi_reg}
    \{e_1, ..., e_{2n}\}:=\{ v_1, ..., v_n, h_1, ..., h_n\}
  \end{equation}
  of $T$ at $v \in TM$ where $v_1, ..., v_n$ are vertical, $h_1, ...,
  h_n$ are horizontal, the last $2n-k$ vectors generate $\H$ and
  $v_k:=C$. Then we have
  \begin{alignat*}{1}
    g_1(\Pfff)&:=\mathrm{Ker}\, (\sigma_1(\Pfff))= \{B_1\in T^* \ | \ 
    B_1(e_i)=0, \ i=k,...,2n \},
    \\
    g_2(\Pfff)&:=\mathrm{Ker}\, (\sigma_2(\Pfff))= \{B_2\in S^2T^* \ | \ 
    B_2(e_i, e_j)=0, \ i=1,...,2n, \ j=k,...,2n \}
  \end{alignat*}
  and for $1 \leq m < k$,
  \begin{alignat*}{1}
    g_{1,e_1, ..., e_m}(\Pfff) :&= \mathrm{Ker}\, (\sigma_1(\Pfff))
    \cap \{B_1\in T^* \ | \ B_1(e_i)=0, \ i=1,..,m \}
    \\
    & = \bigl\{B_1\in T^* \ | \ B_1(e_i)=0, \ i\in \{1,..,m\} \cup
    \{k, ..., 2n\} \bigl\}.
  \end{alignat*}
  The dimension of these spaces are
  \begin{alignat*}{1}
    &\mathrm{dim}\,\bigl(g_1(\Pfff)\bigl)= k-1,
    \\
    &\mathrm{dim}\,\bigl(g_2(\Pfff)\bigl)= \frac{k(k-1)}{2},
    \\
    &\mathrm{dim}\,\bigl( g_{1,e_1, ..., e_m}(\Pfff)\bigl) = \left\{
      \begin{aligned}
        k-&1-m, & & \quad \mathrm{for} \quad m=1,...,k-1,
        \\
        &0, & & \quad \mathrm{for} \quad m=k,...,2n,
      \end{aligned}
    \right.
  \end{alignat*}
  therefore
  \begin{alignat*}{1}
    \mathrm{dim}\,\bigl(g_1(\Pfff)\bigl) & + \sum_{m=1}^{2n}
    \mathrm{dim}\,\bigl( g_{1,e_1, ..., e_m}(\Pfff)\bigl) = (k-1) +
    \sum_{m=1}^{k-1} (k-1-m) = \frac{(k-1)k}{2}
    \\
    & = \mathrm{dim}\,\bigl(g_2(\Pfff)\bigl).
  \end{alignat*}
  This shows that (\ref{quasi_reg}) is a quasi-regular basis for
  $\Pfff$. The existence of such a basis proves Lemma \ref{lemma_2}.
\end{proof}

\section{Landsberg metrizability}
\label{cha:landsberg}

In this paragraph we will investigate the following problem:
\textit{under what conditions can a given second order differential
  equation (\ref{eq:sode}) be the geodesic equation of a Finsler
  metric of Landsberg type?}  As we explained in Proposition
\ref{theo:reformulation}, to answer this question we have to look for
a solution of the PDE system consisting of (\ref{eq:E_hom_loc}),
(\ref{eq:E_EL_loc}) and (\ref{eq:E_h_metric_loc}).  Let us consider
the third order system
\begin{equation}
  \label{eq:e_L_hom_g}
  \Pl = (\Pc , \, \Pe, \, \Pg)
\end{equation}
where $\Pc$, $\Pg$ and $\Pe$ are defined by (\ref{eq:E_hom}),
(\ref{eq:E_h_metric_P}) and (\ref{eq:E_EL_P}).  We will prove the
following theorems:

\medskip

\begin{theorem}
  \label{theo:landsberg} 
  (Reduction of $\Pl$) The third-order partial differential system
  $\Pl E=0$ is equivalent to the first order system
  \begin{equation}
    \label{eq:L_syst_2_1}
    \left\{\quad
      \begin{aligned}
        {\mathcal L}_CE-2E&=0,
        \\
        d_{\l} E & =0,
        \end{aligned}
    \right.
  \end{equation}
  where $\mathfrak L$ is the distribution generated by the horizontal
  vector fields, the image of the Berwald curvature and their
  successive Lie-brackets and $\l : TTM \to \L$ is an arbitrary
  projection of $TTM$ onto $\L$.
\end{theorem}

\begin{remark_num}
  The second equation of (\ref{eq:L_syst_2_1}) means simply that the
  Lie-derivative of $E$ with respect to vector-fields in the
  distribution $\L= \mathrm{Im}\, \l$ is zero.  This property is
  independent of the projection $\l$ of $\L$ chosen.
\end{remark_num}  

\begin{remark_num}
  $E$ is a solution of (\ref{eq:L_syst_2_1}) if and only if it is a
  solution of
    \begin{equation}
      \tag{\ref{eq:L_syst_2_1}'}
      \label{eq:L_syst_2_1_1}
      \left\{\quad
        \begin{aligned}
          {\mathcal L}_CE-2E&=0,
          \\
          d_h E & =0,
          \\
          d_{\mathcal R} E & =0,
        \end{aligned}
      \right.
    \end{equation}
    where $h$ is simply the horizontal projection associated to
    $\Gamma$. However, under the assumption of regularity, the system
    (\ref{eq:L_syst_2_1}) is integrable but (\ref{eq:L_syst_2_1_1}) in
    general is not, because it is not containing its compatibility
    conditions.
  \end{remark_num}

\begin{theorem}
  \label{theo:landsberg_2}
  Let $S$ and $M$ be analytical, and suppose that rank of $\L$
  constant in a neighborhood of $v \in \TM$. Then there exists a
  Finsler metric of Landsberg type in a neighborhood of $v$ whose
  geodesics are given by $S$, if and only if, $C \not \in \L$, and the
  kernel of the first prolongation of (\ref{eq:L_syst_2_1}) at $v$
  contains a positive definite initial condition.
\end{theorem}

In order to prove the above theorems, we need the following

\begin{lemma}
  \label{theo:lemma_3}
  Let us consider the differential operator $d_\P: C^\infty(TM) \to
  {\mathcal S}ec\,(S^3T^*)$, where $\P$ is the Berwald curvature.  For all
  $X,Y,Z \in T$ we have
  \begin{equation}
    \label{eq:3}
    \Pg E \, (X,Y,Z) = \n^2 \PH E \, (JY, JZ,hX) + d_\P E \, (X, Y, Z),
  \end{equation}
  where $\nabla$ is the Berwald connection and $\PH$ is introduced in
  (\ref{eq:d_p}).
\end{lemma}

\textit{Proof.} The three terms in (\ref{eq:3}) are all
semi-basic in $X$, $Y$ and $Z$.  Putting $X = \p{x^i}$, $Y = \p{x^j}$
and $Z = \p {x^k}$, we have
\begin{alignat*}{1}
  & (\Pg E) \Bigl( \p {x^i}, \p {x^j}, \p {x^k}\Bigl)- (\n^2 \PH E)
  \Bigl(\p {y^j}, \p {y^k}, h \Bigl(\p {x^k}\Bigl)\Bigl)
  = \frac{\partial ^3 E}{\partial x^i \partial y^j \partial y^k}
  \\
  & \quad \phantom{=} - \Gamma_i^l \frac{\partial ^3 E}{\partial
    y^l \partial y^j \partial y^k} - \Gamma_{ij}^l
  \frac{\partial^2 E}{\partial y^l \partial y^k} -
  \Gamma_{ik}^l \frac{\partial^2 E}{\partial y^l \partial
    y^j} - \p {y^j} \p {y^k} \left( \frac{\partial E}{\partial x^i } -
    \Gamma_i^l \frac{\partial E}{\partial y^l} \right)
  \\
  & \quad = \frac{\partial ^2 \Gamma^l_{i}}{\partial y^j \partial
    y^k} \frac{\partial E}{\partial y^l}
  = -\frac{1}{2}\frac{\partial^3 f^l}{\partial y^i \partial y^j
    \partial y^k} \frac{\partial E}{\partial y^l}
  = d_\P E \, \Bigl(\p {x^i}, \p {x^j}, \p {x^k}\Bigl).
\end{alignat*}
\newline \hphantom{.} \hfill{$\Box$}

\noindent
\begin{idfn}
  Since $\h \circ h = h$, we have $\PH(hX) = d_hE(hX)$, and we have
  the relation
  \begin{equation}
    \tag{\ref{eq:3}'} \Pg E \, (X,Y,Z) = \n^2 d_hE \, (JY, JZ,hX) +
    d_\P E \, (X, Y, Z).
  \end{equation}
  expressed in terms of the horizontal projection $h$.
\end{idfn} 

\smallskip

\noindent
\textit{Proof of Theorem \ref{theo:landsberg}}.

1) If $E:TM \to \R$ is a solution of $\Pl=(\Pc, \, \Pe, \, \Pg)$, then
by Theorem \ref{theo:finsler} we obtain that $\PH E =0$.  In
particular, $E$ is constant on the horizontal distribution.  Moreover,
we can find from (\ref{eq:3}) that $d_\P E=0$, i.e.~$E$ is constant on
the image of the Berwald curvature.  Consequently $E$ has to be
constant on the distribution $\L$ generated by the horizontal
vector-field and the image of Berwald curvature.
  
2) Conversely, let $E:TM \to \R$ be a solution of
(\ref{eq:L_syst_2_1}). By the construction $\H \subset \L$, we obtain
that $E$ is a solution of $\PH$. By Theorem \ref{theo:finsler}, $E$ is
also a solution of $\Pc$ and $\Pe$.  Moreover, $\mathrm{Im} \, \P
\subset \L$ implies $d_\P E=0$ and by Lemma \ref{theo:lemma_3} we have
$\Pg E=0$.  Therefore $E$ is a solution of the system $\Pl=(\Pc, \,
\Pe, \, \Pg)$.

\smallskip

By 1) and 2) we conclude that the system $\Pl=(\Pc, \, \Pe, \, \Pg)$
is equivalent to the $1^{\mathrm{st}}$ order system
(\ref{eq:L_syst_2_1}) which proves Theorem \ref{theo:landsberg}.
\hfill{$\Box$}

\medskip

\nn \textit{Proof of Theorem \ref{theo:landsberg_2}.}

The reasoning is completely analogous to the proof of Theorem
\ref{theo:finsler_2}. Indeed, $E:TM \to \R$ is a Landsberg-type
Finsler metric associated to $S$, if and only if, $g_E$ is positive
definite and $E$ is a solution of the system $\Pl =(\Pc, \, \Pe , \,
\Pg )$.

If $C \in \L$ and $E:TM \to \R$ is a solution of
(\ref{eq:L_syst_2_1}), then $dE=0$. So $E$ is not a regular
Lagrangian, and $S$ cannot be variational.

Suppose that $C \not \in \L$ and that $\L$ has constant rank in a
neighbourhood of $v\in \TM$.  By Theorem \ref{theo:landsberg} we know
that $E$ is a solution of $\Pl$ if and only if it is a solution of
(\ref{eq:L_syst_2_1}).  Therefore, it is sufficient to consider this
first order PDE and show that it has a solution.

By the hypotheses, $\L$ is of constant rank in a neighbourhood of
$v\in TM$, the system (\ref{eq:L_syst_2_1}) is regular.

A computation, completely analogous to that of made in the proof of
Proposition \ref{theo:D_int}, shows that (\ref{eq:L_syst_2_1}) is
formally integrable.  Consequently, for every initial condition, there
exists an analytical solution to (\ref{eq:L_syst_2_1}) in a
neighbourhood of $v \in \TM$. Using Theorem \ref{theo:landsberg}, this
function will be a solution of the system $\Pl=(\Pc, \, \Pe , \, \Pg
)$, and therefore it will be a Landsberg type Finsler metric in a
neighborhood of $v$ with geodesics determined by $S$.
\hfill{$\Box$}

\section{Remarks and examples of Finsler and Landsberg metrizability}
\label{sec:remarks}

Theorems \ref{theo:finsler}, \ref{theo:C_D}, \ref{theo:finsler_2},
\ref{theo:landsberg} and \ref{theo:landsberg_2} give us a powerful
method to test the metrizability of a second order ordinary
differential system. We mention here only some direct consequences.

\begin{proposition}
  \label{cor:remarks}
  A quadratic second order differential equation is Landsberg
  metrizable if and only if it is Finsler metrizable.
\end{proposition}

Indeed, in the quadratic case, the functions $f^i(x, \dot{x})$ are
quadratic in the $\dot{x}$ variable and the Berwald curvature $\P$
vanishes identically.  Therefore the distribution $\L$ coincides with
$\H$.  \hfill{$\Box$}

\begin{theorem}
  \label{theo:L=T}
  If $\mathrm{rank} \, \L=2n$ (resp.~$\mathrm{rank}\,\H=2n$), then the
  spray is not Landsberg (resp.~Finsler) metrizable.
\end{theorem}

Indeed, in this cases $\L=T$ (resp.~$\H=T$). If $E:TM \to \R$ is a solution
of (\ref{eq:L_syst_2_1}) (resp. (\ref{eq:syst_2_1})), then $dE=0$, and $E$
cannot be a regular Lagrangian. 
\hfill{$\Box$}

\medskip

\nn \textsl{Examples}
\begin{enumerate}
\item For a generic spray, the image of the curvature $R$ and the image of
  $\P$ generate the whole vertical space. In this case $\L=T$, and
  therefore there is no a regular solution to (\ref{eq:L_syst_2_1}). 
  
\medskip

\item In some cases, even if the image of the curvature $R$ and the
  image of $\P$ do not generate the whole vertical space, nevertheless
  $\L=T$. For example let $f(t):=a \sqrt{t^2+bt+c}$ with $a$, $b$, $c$
  nonzero reals, and consider the system
  \begin{equation}
    \label{ex:3}
    \ddot{x}_1={\dot x}_1^2 \, f \left( \frac{{\dot x}_2}{{\dot x}_1}
    \right) , \quad \ddot{x}_2={\dot x}_1{\dot x}_2 \, f \left(
      \frac{{\dot x}_2}{{\dot x}_1} \right).
  \end{equation}
  In this case
  \begin{math}
    \mathrm{Im} \, \P = \mathrm{Im} \, R
  \end{math}
  is a 1-dimensional distribution of $T$. However, by computing the
  Lie-brackets of horizontal vector fields with the generator of
  $\textrm{Im}\, \P$ we find that
  \begin{displaymath}
    \Bigl[h \p{x^i}, \P \Bigl(\p {x^1}, \p {x^1}, \p
    {x^1}\Bigl)\Bigl]_{i=1,2} \in \textrm{Im} \, \P \quad
    \Leftrightarrow \quad b=c=0.
  \end{displaymath}
  Therefore we have $\mathrm{dim}\, \L =4$, so $\L =T$ and
  (\ref{ex:3}) is not Landsberg metrizable.
\end{enumerate}

\begin{theorem}
  \label{theo:constant} 
  If $\H$ (resp.~$ \L$) contains the vertical lift of a non-zero
  vector field on $M$, then the spray is not Finsler (resp.~Landsberg)
  metrizable.
\end{theorem}

\begin{proof}
  Let $X \in \H$ (resp.~$X \in \L$) be a vertical lift, namely
  $X=Z^v$, $Z \in \mathfrak{X}(M)$. Then, locally, $X = (X^\alpha
  \circ \pi) \frac{\partial}{\partial y^\alpha}$, where the functions
  $X^\alpha$ are defined on a domain of $M$.  If $S$ is Finsler
  (resp.~Landsberg) metrizable, then the corresponding energy function
  $E:TM \to \R$ is a regular Lagrangian, and it is a solution of $\Pf$
  (resp.~$\Pl$). Using Theorem \ref{theo:finsler} (resp.~Theorem
  \ref{theo:landsberg}) we get that $E$ is a solution of
  (\ref{eq:syst_2_1}) (resp.~(\ref{eq:L_syst_2_1})), and in particular
  $E$ is constant on every vector fields of $\H$ (resp.~$\L$).
  
  Since $X \in \H$ (resp.~$X \in \L$) we have ${\mathcal L}_XE=0$.
  Taking the derivatives with respect to the vertical directions and
  using the special form of $X$ we obtain that
  \begin{displaymath}
    0=\frac{\partial {\mathcal L}_X E}{\partial {y^i}}
    = \p {y^i} \left( (X^j \circ \pi) \frac{\partial E}{\partial y^j }
    \right)
    = (X^j \circ \pi) \frac{\partial^2 E}{\partial y^j \partial
      {y^i}},
  \end{displaymath}
  so $E$ cannot be a regular Lagrangian. This contradicts the
  hypothesis. Therefore the spray is not Finsler (resp. Landsberg)
  metrizable.
\end{proof}

\medskip

\nn \textsl{Example.}  Let us consider the system
\begin{equation}
  \label{eq:ex_1}
  \left\{
    \begin{aligned}
      \ddot{x}^1:= \lambda_1(x) \, f(x, \dot{x}),
      \\
      \ddot{x}^2:= \lambda_2 (x)\, f(x, \dot{x}),
    \end{aligned}
  \right.
\end{equation}
where $f(x, y)$ is an arbitrary second order homogeneous but
non-quadratic function in $y=(y^1, y^2)$ and $\lambda_1$, $\lambda_2$
arbitrary functions of $x=(x^1, x^2)$.  In this case the image of the
Berwald curvature is generated by the vertically lifted vector field
$X= \lambda_1 \p {y^1} + \lambda_2 \p {y^2}$.  Thus, by Theorem
\ref{theo:constant}, the system is not Landsberg metrizable.

\section{On the existence of non-Berwald type Landsberg spaces}
\label{sec:exist-non-berw}

A Landsberg metric is said to be of Berwald type if the connection
$\Gamma$ is linear, that is in its geodesic equations $\ddot{x}^i =
f^i(x, \dot{x})$ the functions are quadratic in $\dot{x}$.  These
types of spaces can be characterized in terms of the Berwald
curvature: a Landsberg space is of Berwald type if and only if the
Berwald curvature introduced in (\ref{eq:P}) vanishes. One of the most
exciting questions in Finsler geometry is the following:
\begin{center}
  \textit{Are there any non-Berwald Landsberg metrics on a manifold?}
\end{center}
To answer this question a promising strategy is to investigate the
solvability of the system $\Pl=0$, which is a third order differential
system.  Theorems \ref{theo:landsberg} and \ref{theo:landsberg_2} can
be useful for this purpose, because they provide a reduction of $\Pl$
to a much simpler first order differential system.  Far from exploring
fully the possibilities offered by the above theorems, we shall be
content to make the following observations.

\begin{theorem}
  There is no a nontrivial analytic function $f$ such that the
  equations
  \begin{equation}
    \begin{aligned}
      \ddot{x}_1 & = {\dot x}_1^2 \, f({\dot x}_2/{\dot x}_1),
      \\
      \ddot{x}_2 & = {\dot x}_1{\dot x}_2 \, f({\dot x}_2/{\dot x}_1)
    \end{aligned}
  \end{equation}
  constitute the geodesic system of a non-Berwald type Landsberg
  metric.
\end{theorem}

\begin{proof}
  If $f \not \equiv 0$, then $\P \neq 0$.  Unless $f$ satisfies the
  equation
  \begin{math}
    3 f'' f' + f f''' =0,
  \end{math}
  we have   
  \begin{math}
    \mathrm{Im} \, \P \neq \mathrm{Im} \, R.
  \end{math}
  In this case $\L$ is the entire second tangent bundle $TTM$, and
  consequently there is no corresponding Landsberg metric.
    
  If $f$ satisfies the above equation, then it has the form
  \begin{math}
    f(t) = a \sqrt{t^2+bt+c}
  \end{math}
  with $a,b,c \in \R$.  Computing the Lie brackets $[h (\p {x^i}),
  \P(\p {x^1}, \p {x^1}, \p {x^1})]$ we find that they are in the
  subspace generated by $\mathrm{Im} \, \P$ if and only if $b=c=0$.
  But in this case $\P=0$ which contradicts our hypotheses.
\end{proof}

\medskip

\begin{proposition}
  The system
  \begin{equation}
    \label{ex:4}
    \begin{aligned}
      \phantom{\qquad a \in \R, \ t \in \Bbb N,} \ddot{x}_1&=a \,
      {\dot x}_1^{2-t}\, {\dot x}_2^t, \qquad a \in \R, \ t \in \Bbb N,
      \\
      \phantom{\qquad b \in \R, \ s \in \Bbb N,} \ddot{x}_2&=b \,
      {\dot x}_1^{2-s} \, {\dot x}_2^s , \qquad b \in \R, \ s \in \Bbb N,
    \end{aligned}
  \end{equation}
  cannot be the geodesic system of a non-Berwald type Landsberg
  metric.
\end{proposition}

\begin{proof}
  Let us consider the spray $S$ corresponding to (\ref{ex:4}):
\begin{displaymath}
  S = y^1 \p {x^1}+ y^2 \p {x^2} + a \, y_1^{2-t}\, y_2^t \p {y^1}+b
  \, y_1^{2-s}\, y_2^s \p {y^2}.
\end{displaymath}
If $s \in \{0,1,2\}$ or $t \in \{0,1,2\}$, then $S$ is not the
geodesic equation of a non-Berwald type Landsberg metric.

\smallskip

Indeed, in case of $s, t \in \{0,1,2\}$, then $\P=0$ and therefore the
Berwald connection is linear. If the system is Finsler-metrizable,
then it is also Landsberg metrizable (Corollary \ref{cor:remarks}),
and the corresponding Finsler space is of Berwald type.  If $s \in
\{0,1,2\}$ or $t \in \{0,1,2\}$, then $\mathrm{Im} \, \P$ generated by
$\p {y^1}$ or $\p {y^2}$.  As we explained in Theorem
\ref{theo:constant}, in these cases there is no regular Lagrangian
associated to the system.
  
\medskip

If $s, t \not \in \{0,1,2\}$, then the image of $\P$ is generated by
the vector-fields
\begin{displaymath} 
  \p {y_1} + \frac{b s (s-1)(s-2)y_2^{s-t} }{ a t (t-1)(t-2)y_1^{s-t}}
  \p {y_2}.
\end{displaymath}
If in addition $s= t$, then using Theorem \ref{theo:constant} we
obtain that there is no regular Lagrangian associated to the system.
Let us suppose now that $s \neq t$. The image of the curvature $R$ is
generated by the vector field
  \begin{alignat*}{1}
    \p {y_1} + \frac{b}{a} \frac{
      y_2^{s} y_{1}^{2 t-s+2}\, b s(2-s)
      + y_2^{t+1} y_{1}^{t+1} \, a(2s^2-4s+2t-st)
    } {
      y_2^{t+1} y_{1}^{t+1} \, a t(t-2)
      + y_2^s y_{1}^{2t-s+2} \, b (st-2t^2 + 2t) }
    \left( \frac{y_2}{y_1} \right)^{s-t} \p {y_2}.
\end{alignat*}
If $s \neq t+1$, or $s \neq t-1$, then $\mathrm{Im} \, R \neq
\mathrm{Im} \, \P$. Since $ T^h \oplus \mathrm{Im} \, R \oplus
\mathrm{Im} \, \P \subset \L$ we obtain that $\L=T$. Using Theorem
\ref{theo:L=T} we find that $S$ is not Landsberg-metrizable.

If $s = t+1$ or $s = t-1$, then $\mathrm{Im} \, R = \mathrm{Im} \,
\P$. Computing the Lie-brackets of the horizontal vector-fields with
the image of the Berwald curvature we find that
\begin{displaymath}
  \bigl[ T^h, \, \mathrm{Im} \, \P \bigl] \, \nsubseteq \, T^h \oplus
  \mathrm{Im} \, \P.
\end{displaymath}
We arrive at $\L=T$, and using Theorem \ref{theo:L=T} we conclude
again that $S$ is not Landsberg-metrizable.  
\end{proof}

\end{document}